\documentclass{article}
\usepackage{mathrsfs}[12pt]
\usepackage{bm, amscd,  amsfonts, amssymb, cite,texdraw}
\usepackage{amsmath, epsfig, cite, color, colordvi}
\usepackage{latexsym,graphicx}

\makeatletter \@addtoreset{equation}{section} \makeatother

\newtheorem{theo}{Theorem}[section]

\newtheorem{coro}[theo]{Corollary}

\def\0{\mathbf 0}

\def\mid{{\,|\,}}

\def\pf{\noindent {\it Proof.\ }}
\def\qed{\hfill \rule{4pt}{7pt}}
\begin{document}

\begin{center}
{\Large    Combinatorial Telescoping for an Identity of
          Andrews\\[2mm]

          on Parity in Partitions}

\vskip 3mm

 William Y.C. Chen$^1$,
Daniel K. Du$^2$ and Charles B. Mei$^3$\\[5pt]
Center for Combinatorics, LPMC-TJKLC \\
Nankai University, Tianjin 300071, P. R. China

\vskip 3mm

 E-mail: $^1$chen@nankai.edu.cn,
$^2$dukang@mail.nankai.edu.cn, $^3$meib@mail.nankai.edu.cn

\end{center}

\begin{abstract}
Following the method of combinatorial telescoping for alternating sums
 given by Chen, Hou and
Mu, we present a combinatorial telescoping approach to partition
identities on sums of positive terms. By giving a classification of
the combinatorial objects corresponding to a sum  of positive terms,
we establish bijections that lead
 a telescoping relation. We illustrate this idea by
giving a combinatorial telescoping relation for a classical
identity of MacMahon. Recently, Andrews posed a problem of finding a
combinatorial proof of an identity on the $q$-little Jacobi
polynomials which was derived based on a recurrence relation.
We find a combinatorial classification of certain triples of
partitions and a sequence of bijections. By the method of cancelation,
we see that there exists
 an involution for a recurrence relation that implies the identity of Andrews.

\end{abstract}



\section{Introduction}

In his study of parities in partition identities, Andrews \cite{Andrews2010}
obtained the following identity  on the little $q$-Jacobi
polynomials \cite[p. $27$]{Gasper}:
\begin{equation}\label{que-12}
 _2\phi_1\left(\!
\begin{array}{c}
  q^{-n}, q^{n+1}\\
  -q
\end{array}\!;q,-q
\right)
=(-1)^nq^{{n+1 \choose
2}}\sum_{j=-n}^{n}(-1)^jq^{-j^2}.
\end{equation}
Let $G_{n}(q)$ denote the sum on the left hand side of (\ref{que-12}).
 Andrews \cite{Andrews2010} established the following recurrence
relation for $n\geq 1$,
\begin{equation}\label{gn}
G_n(q) +q^n G_{n-1}(q) = 2q^{-{n \choose 2}},
\end{equation}
from which \eqref{que-12} can be easily deduced. As one of the fifteen open problems, Andrews
asked for a combinatorial proof of the above identity (\ref{que-12}).

In this paper, we give a combinatorial treatment of a
 homogeneous recurrence relation for the sum
\[ F_n(q)=q^{n \choose 2}\,_2\phi_1\left(\!
\begin{array}{c}
  q^{-n}, q^{n+1}\\
  -q
\end{array}\!;q,-q
\right),
\]
which is a consequence of  recurrence relation (\ref{gn}).
More precisely, for $n\geq 2$ we have
\begin{equation}\label{rec-fn}
F_n(q)+(q^{2n-1}-1)F_{n-1}(q)-q^{2n-3}F_{n-2}(q)=0.
\end{equation}
It is readily seen that identity \eqref{que-12}
is an immediate consequence of  (\ref{rec-fn}).

The main objective of this paper is to present a combinatorial treatment
of the recurrence relation (\ref{rec-fn}).
To this end, we present the method of combinatorial telescoping for
sums of positive terms, which is a variant of the method of combinatorial
telescoping for alternating sums. In this framework, we find a
classification of certain triples of partitions and a sequence of bijections, leading to
a combinatorial proof of the above recurrence relation \eqref{rec-fn}.

Recall that Chen, Hou and Mu \cite{Chen2005} presented the method of
combinatorial telescoping for alternating sums. Consider the
alternating sum
\begin{equation}\label{eq-ct11}
  \sum_{k=0}^\infty (-1)^k f(k).
\end{equation}
A combinatorial telescoping   for the
above alternating sum means a  classification of certain
combinatorial objects along with a sequence of bijections. This
method can be used to show that the above alternating sum satisfies
a recurrence relation, and it applies to many  $q$-series identities on
alternating sums such as Watson's
identity \cite{Watson}
\begin{equation}\label{Watson}
\sum^{\infty}_{k=0} (-1)^k \frac{1-aq^{2k}}{(q;q)_k(aq^{k};q)_{\infty}}a^{2k}q^{k(5k-1)/2}=
\sum^{\infty}_{n=0}\frac{a^nq^{n^2}}{(q;q)_n},
\end{equation}
and Sylvester's identity \cite{Syl}
\begin{equation}\label{Syl}
\sum^{\infty}_{k=0} (-1)^k q^{k(3k+1)/2}x^k \frac{1-xq^{2k+1}}{(q;q)_k(xq^{k+1};q)_{\infty}}=1.
\end{equation}

In this paper, we consider a summation of the following form
\begin{equation}\label{sum}
 \sum_{k=0}^\infty f(k).
\end{equation}
Suppose that $f(k)$ is a weighted count of a set $A_k$, that is,
\[
 f(k) = \sum_{\alpha\in A_k} w(\alpha).
\]
We wish  to find sets $B_k$ and $H_k$ $(k\geq 0)$ with a weight assignment $w$ such
that there exists a weight preserving bijection
\begin{equation}\label{biject-ct2}
 \phi_k \colon A_k \cup H_k \rightarrow B_k  \cup H_{k+1},
\end{equation}
where $\cup$ means disjoint union. Let $g(k)$ and
$h(k)$ be the weighted count of the sets $B_k$ and $H_k$,
respectively, namely,
\[
 g(k) = \sum_{\alpha\in B_k} w(\alpha) \quad \mbox{and} \quad
 h(k) = \sum_{\alpha\in H_k} w(\alpha),
\]
then the bijection \eqref{biject-ct2} implies the relation
\begin{equation} \label{rec-ct2}
 f(k)+ h(k)  = g(k) + h(k+1).
\end{equation}
Just like the conditions for   creative telescoping \cite{Zeilberger1991,Knuth1994,A=B1996}, we suppose
that $H_0=\emptyset$ and $H_k$ vanishes for sufficiently
large $k$. Summing \eqref{rec-ct2} over
$k$ gives the following relation
\begin{equation}\label{eq-ct21}
\sum_{k=0}^\infty  f(k) =  \sum_{k=0}^\infty  g(k),
\end{equation}
which is equivalent to a recurrence relation of the sum (\ref{sum}).

Combining all the bijections $\phi_k$ in \eqref{biject-ct2}, we get a correspondence
\begin{equation}\label{biject-ct-main}
 \phi \colon  A\cup H \longrightarrow B \cup H,
\end{equation}
given by $\phi(\alpha)=\phi_k(\alpha)$ if $\alpha\in A_k\cup H_k$,
where
\[
 A =\bigcup_{k=0}^\infty A_k,\quad
 B =\bigcup_{k=0}^\infty B_k,\quad
 H =\bigcup_{k=0}^\infty H_k.
\]

 By the method of cancelation, see Feldman and Propp
 \cite{FP},
the above bijection $\phi$ implies a bijection
\[
\psi\colon  A\longrightarrow B.
 \]
More specifically,  we can define the bijection $\psi\colon A\rightarrow B$ by setting $\psi(a)$
 to be the first element $b$ that falls into $B$ while iterating the action of
 $\phi$ on $a\in A$.

For the purpose of this paper, we shall express $A_k$ as a sum over $n$,
namely,
\[
 A_k= \bigcup_{n=0}^\infty A_{n,k}.
\]
It should be noted that our bijections do  not require an  explicit formula
for $A_{n,k}$. Roughly speaking, our idea is to use bijections to establish a telescoping relation involving $A_{n,k}$ possibly with  coefficients depending only on $n$.

For any $n$ and $k$,  we aim to find bijections $\phi_{n,k}$ for a given integer $r$:
\begin{equation}\label{biject-recurence}
\phi_{n,k} \colon \bigcup_{i=0}^r \{a_i(n)\} \times A_{n-i,k}\cup
H_{n,k} \rightarrow  \bigcup_{i=0}^r \{b_i(n)\} \times A_{n-i,k}\cup
H_{n,k+1},
\end{equation}
where the leading coefficients of $a_i(n)$ and  $b_i(n)$ are
positive, and $ \{0\}\times A_{n-i,k}$ is considered as the empty set. Let
\[
 F_{n,k} = \sum_{\alpha\in A_{n,k}} w(\alpha)
\]
be a weighted count of the set $A_{n,k}$, and let
\[
 F_n = \sum_{k=0}^\infty F_{n,k}.
\]
Indeed, the motivation to find the bijections given in \eqref{biject-recurence}
is to obtain a recurrence
relation of $F_n$. Once the relation \eqref{rec-fn} is established, we immediate
get  \eqref{que-12}.

This paper is organized as follows. In Section \ref{sec-MacMahon},
we illustrate our method of combinatorial telescoping for sums of positive
terms by giving a telescoping proof of an identity of
MacMahon \cite[p.$41$]{Pak2006}. In Section \ref{sec-q12}, we provide a
 solution to Problem $12$ of Andrews \cite{Andrews2010} by using
 the idea of combinatorial telescoping to construct the recurrence
relation \eqref{rec-fn} for the following equivalent form of (\ref{que-12}):
\begin{equation}\label{eq-q12a}
\sum_{k=0}^{n}\frac{(q^{n-k+1};q)_{2k}}{(q^2;q^2)_k}q^{{n-k \choose
2}}=(-1)^nq^{n^2}\sum_{j=-n}^{n}(-1)^jq^{-j^2},
\end{equation}
which can be obtained by multiplying  both sides of
\eqref{que-12} by $q^{n \choose 2}$.

\section{MacMahon's identity}\label{sec-MacMahon}

In this section, we use MacMahon's identity on partitions
 to illustrate the idea of combinatorial
telescoping for sums of positive terms.

Let us recall some notation and definitions  in
\cite{Andrews1998}. A \emph{partition} is a nonincreasing
finite sequence of positive integers
$\lambda=(\lambda_1,\ldots,\lambda_{\ell})$. The integers
$\lambda_i$ are called the parts of $\lambda$. The sum of parts and
the number of parts are denoted by $\mid \lambda \mid =\lambda_1+
\dots +\lambda_{\ell}$ and $\ell(\lambda)=l$, respectively. The
special partition with no parts is denoted by $\emptyset$. Denote by
${D}$ the set of partitions with distinct parts, and denote
by  $
{E}$  the set of partitions with even parts. We shall use
diagrams to represent partitions and use rows to represent parts.

We shall adopt the common notation and terminology on basic hypergeometric series in \cite{Gasper}.
The $q$-shifted factorials and the $q$-binomial coefficients, or the Gaussian coefficients,
are defined by
\[
(a;q)_n=(1-a)(1-aq)\cdots(1-aq^{n-1}),\quad
(a;q)_\infty=\prod_{i=0}^\infty (1-aq^i),
\]
and
\[
 \left[\!
\begin{array}{c}
  n \\
  k
\end{array}\!
\right]_q =\frac{(q;q)_n}{(q;q)_k(q;q)_{n-k}}.
\]

In his classical treatise \cite{MacMahon},
MacMahon gave combinatorial proof of the following identity,
see also Pak \cite[p $41$]{Pak2006}:
\begin{equation}\label{eq-MacMahon}
 \sum_{k=-m}^n z^k q^{k^2}
 \left[\!
\begin{array}{c}
  m+n \\
  m+k
\end{array}\!
\right]_{q^2} = (-q/z,q^2)_m(-zq;q^2)_n.
\end{equation}
It is easily seen that as $m,n \rightarrow \infty$, MacMahon's identity reduces to
Jacobi's triple product identity \cite{Gasper}.

To prove the identity \eqref{eq-MacMahon}, we first give
a combinatorial telescoping argument for the following recurrence
\begin{equation}\label{rec-ex1}
 F_{n,m}(q)=(1+q^{2m-1}/z) F_{n,m-1}(q),
\end{equation}
where $F_{n,m}(q)$ denotes the sum on the left hand side of
\eqref{eq-MacMahon}. To compute $F_{n,m}(q)$, we still need the initial value $F_{n,0}(q)$.

Again, by combinatorial telescoping we get the following recurrence for $F_{n,0}(q)$:
\begin{equation}\label{rec-Mac-initial}
 F_{n,0}(q)=(1+zq^{2n-1}) F_{n-1,0}(q).
\end{equation}

Now we construct bijections for the recurrence relation \eqref{rec-ex1}.
 For a positive integer $k$, we denote the square partition with  $k$ rows by $S_k$,
 namely, the partition with $k$ occurrences of the part $k$. For $k=0$, $S_k$
 is considered as the empty partition. Moreover,
  we define $S_{-k}$  to be the square partition with $k$ rows associated with a minus sign.
  We   call $S_k$ a positive square partition, and
  call $S_{-k}$ a negative square partition.

   To give a combinatorial interpretation of the left hand side of \eqref{eq-MacMahon},
   for $-m\leq k \leq n$, we define the following set
of pairs of partitions
\[
 P_{n,m,k} =
 \left\{
 ( \lambda, \mu )\colon
  \lambda = S_k,\;
  \mu_1\le 2m+2k, \; \ell(\mu)\le n-k,\; \mu\in{E}
\right\}.
\]
In other words, $\lambda=S_k$ is a square partition,
$\mu$ is a partition with at most $n-k$ even parts but no odd parts
such that
   the largest part does not exceed
 $2m+2k$.
 It can be easily verified that the $k$-th summand of the left hand side of \eqref{eq-MacMahon} can
be viewed as a weighted count of $P_{n,m,k}$, that is,
\[
\sum_{(\lambda,\mu)\in P_{n,m,k}}z^{k}q^{\mid \lambda \mid+\mid \mu
\mid}=z^k q^{k^2}
 \left[\!
\begin{array}{c}
  m+n \\
  m+k
\end{array}\!
\right]_{q^2}.
\]
Let
\[
 G_{n,m,k} = \{(\lambda,\mu)\in P_{n,m,k}\colon \mu_1=2m+2k\}.
\]
 By definition, $G_{n,m,k} = \emptyset$ for $k <-m$ or $k \geq n$. For
integers $m,n\ge0$ and $-m\le k \le n$,
 we shall construct a bijection
\[
 \phi_{n,m,k} \colon P_{n,m,k} \cup G_{n,m,k-1} \longrightarrow
 P_{n,m-1,k}\cup\{2m-1\}\times P_{n,m-1,k} \cup G_{n,m,k}.
\]
This bijection can be easily deduced from the following
classification of
\[
P_{n,m,k} \cup G_{n,m,k-1} .
\]
Let $( \lambda, \mu)$ be a pair of partitions in $P_{n,m,k} \cup G_{n,m,k-1}$.
\begin{enumerate}
\item For $(\lambda,\mu)\in P_{n,m,k}$, if $\mu_1=2m+2k$, then
$(\lambda,\mu)\in G_{n,m,k}$. We set
$\phi_{n,m,k}(\lambda,\mu)=(\lambda,\mu)$.

\item For $(\lambda,\mu)\in P_{n,m,k}$, if $\mu_1<2m+2k$, we have $\mu_1\le
2m+2k-2$, which implies that $(\lambda,\mu) \in P_{n,m-1,k}$. We set
$\phi_{n,m,k}(\lambda,\mu)=(\lambda,\mu)$.

\item For $(\lambda,\mu) \in G_{n,m,k-1}$,  $\lambda$ is the square partition $S_{k-1}$, we set $\lambda'=S_k$.
  Removing the first row of $\mu$, we obtain $\mu'$. It is easy to check
 that the resulting pair of partitions $(\lambda',\mu')$ belongs
to $P_{n,m-1,k}$. Set
$\phi_{n,m,k}(\lambda,\mu)=(2m-1,(\lambda',\mu'))$.
\end{enumerate}
 Define the
weight function $w$ on $P_{n,m,k}$ and $(2m-1)\times P_{n,m-1,k}$ as
follows
\begin{eqnarray*}
 w(\lambda,\mu) &=& z^{k}q^{\mid\lambda\mid+\mid\mu\mid},
 \\[5pt]
 w(2m-1,(\lambda,\mu)) &=&\frac{q^{2m-1}}{z} z^{k}q^{\mid\lambda\mid+\mid\mu\mid}.
\end{eqnarray*}
 It can be
verified that $\phi_{n,m,k}$ is a weight preserving bijection. This yields
recurrence relation \eqref{rec-ex1}.

We now turn to the evaluation of the initial value $F_{n,0}(q)$. To prove the identity
\begin{equation}\label{eq-Mac-initial}
 \sum_{k=0}^n z^k q^{k^2}
 \left[\!
\begin{array}{c}
  n \\
  k
\end{array}\!
\right]_{q^2} = (-zq;q^2)_n,
\end{equation}
we consider the  set of pairs of partitions
\[
 Q_{n,k} = \{( \lambda, \mu )\colon \lambda =S_k,\;
\ell(\mu)\le k, \; \mu_1 \le 2n-2k, \; \mu\in {E}\}.
\]
Notice that  the $k$-th summand of the left hand side
of \eqref{eq-Mac-initial} can be viewed as a weighted count of
$Q_{n,k}$, that is,
\[
\sum_{(\lambda,\mu)\in Q_{n,k}}z^{\ell(\lambda)}q^{\mid \lambda
\mid+\mid \mu \mid}.
\]
Let
\[
 H_{n,k} = \{(\lambda,\mu)\in Q_{n,k}\colon \mu_1=2n-2k\}.
\]
 By definition, $H_{n,k} = \emptyset$ for $k = 0$ or $k \ge n$. For
any integers $n,k \ge 0$, we shall construct a bijection
\[
 \psi_{n,k} \colon Q_{n,k} \cup H_{n,k+1} \longrightarrow
 Q_{n-1,k}\cup\{2n-1\}\times Q_{n-1,k} \cup H_{n,k}.
\]
This bijection can be easily deduced from the following
classification of
\[
Q_{n,k} \cup H_{n,k+1} .
\]
Let $( \lambda, \mu)$ be a pair of partitions in $Q_{n,k} \cup H_{n,k+1}$.
\begin{enumerate}
\item For $(\lambda,\mu)\in Q_{n,k}$, if $\mu_1=2n-2k$, then
$(\lambda,\mu)\in H_{n,k}$. We set
$\psi_{n,k}(\lambda,\mu)=(\lambda,\mu)$.

\item For $(\lambda,\mu)\in Q_{n,k}$, if $\mu_1<2n-2k$, we have $\mu_1\le
2n-2k-2$, which implies that $(\lambda,\mu) \in Q_{n-1,k}$. We set
$\psi_{n,k}(\lambda,\mu)=(\lambda,\mu)$.

\item For $(\lambda,\mu) \in H_{n,k+1}$, $\lambda$ is the square partition $S_{k+1}$, we set $\lambda'=S_k$. Removing the first row from $\mu$, we obtain
$\mu'$. Clearly, resulting pair of partitions
$(\lambda',\mu')$ belongs to $Q_{n-1,k}$. Set
$\psi_{n,k}(\lambda,\mu)=(2n-1,(\lambda',\mu'))$.
\end{enumerate}
Define the weight function $w$ on $Q_{n,k}$ and $(2n-1)\times
Q_{n-1,k}$ as follows
\begin{eqnarray*}
 w(\lambda,\mu) &=&
 z^{\ell(\lambda)}q^{\mid\lambda\mid+\mid\mu\mid},
 \\[5pt]
 w(2n-1,(\lambda,\mu)) &=& zq^{2n-1}z^{\ell(\lambda)}q^{\mid\lambda\mid+\mid\mu\mid}.
\end{eqnarray*}
One sees that  $\psi_{n,k}$ is a weight preserving bijection. So we get
 the recurrence relation \eqref{rec-Mac-initial}
\begin{equation*}
 F_{n,0}(q)=(1+zq^{2n-1}) F_{n-1,0}(q),
\end{equation*}
where $F_{n,0}(q)$ denotes the sum on the left hand side of
\eqref{eq-Mac-initial}, with the initial value $F_{0,0}(q)=1$.
Since $F_{0,0}=1$, combining the recurrence relations \eqref{rec-ex1} and
\eqref{rec-Mac-initial}, we arrive at MacMahon's identity \eqref{eq-MacMahon}.

\section{An Open Problem of Andrews}\label{sec-q12}

\begin{figure}[h,t]
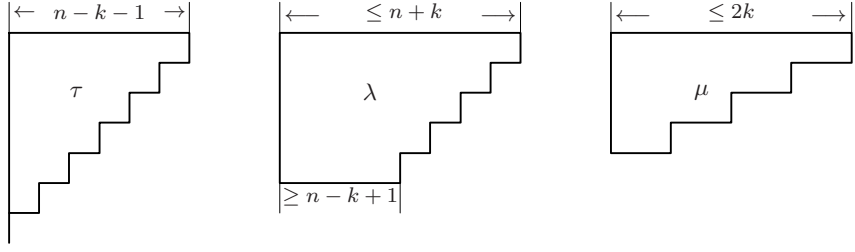

\centertexdraw{ \drawdim mm \linewd 0.25 \setgray 0

\move(0 24)\lvec(24 24)\lvec(24 20)\lvec(20 20 )\lvec(20 16)\lvec(16
16)\lvec(16 12)\lvec(12 12)\lvec(12 8)\lvec(8 8)\lvec(8 4)\lvec(4
4)\lvec(4 0)\lvec(0 0)\lvec(0 24) \move(0 0) \lvec(0 -4)

\textref h:C v:C \htext(12 26.5){\footnotesize$n-k-1$}
 \textref h:C v:C \htext(2 26.5){\footnotesize$\leftarrow$}
  \textref h:C v:C \htext(22 26.5){\footnotesize$\rightarrow$}

\move(36 24)\lvec(68 24)\lvec(68 20)\lvec(64 20)\lvec(64 16)\lvec(60
16)\lvec(60 12)\lvec(56 12)\lvec(56 8)\lvec(52 8)\lvec(52 4)\lvec(44
4)\lvec(36 4)\lvec(36 24)

 \textref h:C v:C \htext(52.5
26.5){\footnotesize$\leq n+k$}
 \textref h:C
v:C \htext(39 26){\footnotesize$\longleftarrow$} \textref h:C v:C
\htext(65 26){\footnotesize$\longrightarrow$}
 \textref h:C v:C
\htext(44 2){\footnotesize$\ge n-k+1$}
 \textref h:C v:C

\move(80 24)\lvec(112 24)\lvec(112 20)\lvec(104 20)\lvec(104
16)\lvec(96 16)\lvec(96 12)\lvec(88 12)\lvec(88 8)\lvec(80
8)\lvec(80 24)

 \textref h:C v:C \htext(96 26.5){\footnotesize$\leq
2k$}
 \textref h:C
v:C \htext(83 26){\footnotesize$\longleftarrow$} \textref h:C v:C
\htext(109 26){\footnotesize$\longrightarrow$}

\linewd 0.1 \setgray 0 \move(36 4)\lvec(36 0) \move(52 4)\lvec(52 0)
\move(36 24)\lvec(36 28) \move(68 24)\lvec(68 28) \move(0 24)\lvec(0
28) \move(24 24)\lvec(24 28) \move(80 24)\lvec(80 28)\move(112
24)\lvec(112 28)

\textref h:C v:C \htext(9 16){\small$\tau$}
 \textref h:C v:C \htext(48 16){\small$\lambda$}
  \textref h:C v:C
\htext(92 16){\small$\mu$} }
  \caption{The diagram $(\tau,\lambda,\mu)\in
P_{n,k}$.}\label{fig-case0}
\end{figure}

In this section, we provide a solution to Problem 12 of
Andrews\cite{Andrews2010} by using the idea of combinatorial
telescoping. Define
\[
P_{n,k}=\left\{
 (\tau,\lambda,\mu)\left | \begin{array}{l}
   \tau=(n-k-1,n-k-2,\ldots,2,1,0),\\[5pt]
   n-k+1 \le \lambda_i \le n+k, \;(i=1,2,\ldots, \ell(\lambda)),\; \lambda \in{D}, \\[5pt]
   \mu_1\leq 2k,  \; \mu \in {E}.
\end{array}
\right.
\right\}.
\]
Figure \ref{fig-case0} gives an illustration of an element of $P_{n,k}$.

In other words, $\tau$ is a triangular partition containing a zero part,
$\lambda$ is a partition with distinct parts, each part of $\lambda$
is smaller than $n+k$ and greater than $n-k+1$, $\mu$ is a partition
with each part even and with the largest part not exceeding $2k$.
As will be seen, we have a reason to include the zero in a triangular partition.
For $k=0$, we have $P_{n,0}=\{(\tau,\emptyset,\emptyset)\}$, where
$\tau=(n-1,n-2,\ldots,2,1,0)$, and for $k>n$, we set $P_{n,k}=\emptyset$.
For  $k=n-1$ and $k=n$, we have
\begin{eqnarray*}
P_{n,n-1}&=&\{(\tau,\lambda,\mu)\colon
   \tau=(0),\;
   2 \le \lambda_i \le 2n-1,
   \; \lambda \in{D},
   \;\mu_1\leq 2n-2,  \; \mu \in {E}\},\\[5pt]
P_{n,n}&=&\{(\tau,\lambda,\mu)\colon
   \tau=\emptyset,\;
   1 \le \lambda_i \le 2n, \; \lambda \in{D},
   \;\mu_1\leq 2n,  \; \mu \in {E}\}.
\end{eqnarray*}
 Notice that we have imposed the distinction between the partition of zero and
 the empty partition. Under this convention,
  one sees that $\bigcup_{k \geq 0}P_{n,k}$ is a disjoint union of $P_{n,k}$.
  Moreover,  the $k$-th
summand $F_{n,k}$ of the left hand side of \eqref{eq-q12a} can be viewed as a
weighted count of $P_{n,k}$, that is,
\[
F_{n,k}=\sum_{(\tau,\lambda,\mu)\in P_{n,k}}(-1)^{\ell(\lambda)}q^{\mid \tau
\mid +\mid \lambda \mid+\mid \mu \mid}.
\]
Notice that the summand term $F_{n,k}$ does not contain the factor $(-1)^k$ as in
an alternating sum. So the summation \eqref{eq-q12a} should be viewed as a sum of positive terms.

Now we give a  combinatorial telescoping relation for
$P_{n,k}$.

\begin{theo}\label{thm-bijection}
For any nonnegative integer $n$ and $0\le k\le n-2$, there is a
bijection
\begin{equation}\label{biject-q12}
 \phi_{n,k}\colon  P_{n,k} \cup \{2n-1\}\times P_{n-1,k-1}
              \rightarrow P_{n-1,k-1} \cup \{2n-3\}\times P_{n-2,k}.
\end{equation}
\end{theo}

\pf For $k=0$, as $P_{n-1,k-1}$ is the empty
set, the bijection $\phi_{n,0}$ is defined by
\[
\phi_{n,0} \colon (\tau,\emptyset,\emptyset) \mapsto
      (2n-3,(\tau',\emptyset,\emptyset)),
\]
where $\tau'$ is obtained from $\tau$ by removing the first two
parts. For example, when $n=2$, $\tau=(1,0)$ and the triple of
partitions is mapped to $(1,(\emptyset,\emptyset,\emptyset))$ belonging to the
set $\{2n-3\}\times P_{n-2,k}$.
Because of the zero part, it is always possible to remove two
parts of $\tau$.

For positive integer $k$, the bijection $\phi_{n,k}$ is essentially
a classification of the set $P_{n,k}$ into four classes, namely,
\[
    P_{n,k}  = A_{n,k}\cup B_{n,k}\cup C_{n,k}\cup P_{n-1,k-1},
\]
where
\begin{align*}
&A_{n,k}=\{(\tau, \lambda, \mu)\in P_{n,k}\colon \lambda_1\le n+k-2,\; \mu_1=2k \}, \\[5pt]
&B_{n,k}=\{(\tau, \lambda, \mu)\in P_{n,k}\colon
\mbox{\small either $n+k$ or $n+k-1$ appears in $\lambda$, but not both} \}, \\[5pt]
&C_{n,k}=\{(\tau, \lambda, \mu)\in P_{n,k}\colon \lambda_1=n+k, \; \lambda_2=n+k-1\}.
\end{align*}
We also need the following classification
\[
  P_{n-2,k}= A'_{n,k}\cup B'_{n,k}\cup C'_{n,k}\cup D_{n,k},
\]
where
\begin{align*}
&A'_{n,k}= \{ (\tau,\lambda,\mu)\in P_{n-2,k}\colon
\lambda_{\ell}\ge n-k+1 \},\\[5pt]
&B'_{n,k}= \{ (\tau,\lambda,\mu)\in P_{n-2,k}\colon
\mbox{\small either $n-k$ or $n-k-1$ appears in $\lambda$, but not both} \},\\[5pt]
 &C'_{n,k}=
\{ (\tau,\lambda,\mu)\in P_{n-2,k}\colon
\lambda_{\ell}= n-k-1 ,\; \lambda_{\ell-1}= n-k ,\; \mu_1=2k \},\\[5pt]
&D_{n,k}= \{ (\tau,\lambda,\mu)\in P_{n-2,k}\colon \lambda_{\ell}=
n-k-1 ,\; \lambda_{\ell-1}= n-k ,\; \mu_1<2k\}.
\end{align*}

Now we are ready to describe the  bijection $\phi_{n,k}$. Assume that $(\tau,\lambda,\mu)$ is a
triple of partitions in $P_{n,k}$.\\

\noindent Case $1$: $(\tau,\lambda,\mu)\in P_{n-1,k-1}$. Set
$\phi_{n,k}(\tau,\lambda,\mu)$ to be $(\tau,\lambda,\mu)$ itself.\\

\begin{figure}[h]
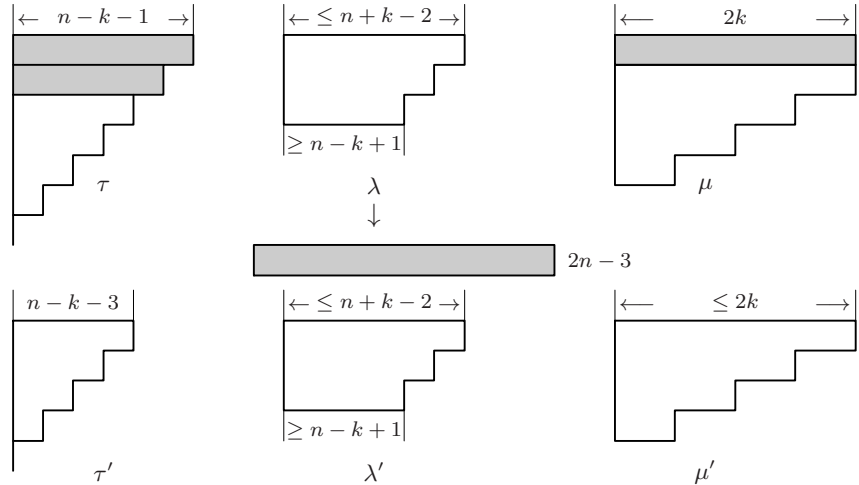

\centertexdraw{ \drawdim mm \linewd 0.25 \setgray 0

\move(0 54) \lvec(20 54) \lvec(20 58)\lvec(24 58)\lvec(24 62)\lvec(0
62)\lvec(0 54) \ifill f:0.8

\move(0 54) \lvec(20 54) \lvec(20 58)\lvec(24 58)\lvec(24 62)\lvec(0
62)\lvec(0 54) \move(0 58)\lvec(20 58) \move(0 54)\lvec(0 38)\lvec(4
38)\lvec(4 42)\lvec(8 42)\lvec(8 46)\lvec(12 46)\lvec(12 50)\lvec(16
50)\lvec(16 54)\move(0 38) \lvec(0 34)

\textref h:C v:C \htext(12 64.5){\footnotesize$n-k-1$}
 \textref h:C
v:C \htext(2 64){\footnotesize$\leftarrow$} \textref h:C v:C
\htext(22 64){\footnotesize$\rightarrow$}

\move(36 62)\lvec(60 62)\lvec(60 58)\lvec(56 58)\lvec(56 54)\lvec(52
54)\lvec(52 50)\lvec(36 50)\lvec(36 62)

 \textref h:C v:C \htext(38 64){\footnotesize$\leftarrow$}
 \textref h:C v:C \htext(48 64.5){\footnotesize$\leq n+k-2$}
 \textref h:C v:C \htext(58 64){\footnotesize$\rightarrow$}

 \textref h:C v:C \htext(44 47.5){\footnotesize$\ge n-k+1$}

\move(80 62)\lvec(112 62)\lvec(112 58)\lvec(80 58)\lvec(80 62)\ifill
f:0.8 \move(80 62)\lvec(112 62)\lvec(112 58)\lvec(80 58)\lvec(80 62)

 \move(80
58)\lvec(112 58)\lvec(112 54)\lvec(104 54)\lvec(104 50)\lvec(96
50)\lvec(96 46)\lvec(88 46)\lvec(88 42)\lvec(80 42)\lvec(80 58)

 \textref h:C v:C \htext(96 64.5){\footnotesize$
2k$}
 \textref h:C
v:C \htext(83 64){\footnotesize$\longleftarrow$} \textref h:C v:C
\htext(109 64){\footnotesize$\longrightarrow$}

\move(32 30)\lvec(72 30)\lvec(72 34)\lvec(32 34)\lvec(32 30) \ifill
f:0.8 \textref h:C v:C \htext(78 32){\footnotesize$2n-3$}

\move(32 30)\lvec(72 30)\lvec(72 34)\lvec(32 34)\lvec(32 30)

\move(0 24)\lvec(16 24)\lvec(16 20)\lvec(12 20 )\lvec(12 16)\lvec(8
16)\lvec(8 12)\lvec(4 12)\lvec(4 8)\lvec(0 8)\lvec(0 24)\move(0 8)
\lvec(0 4)

\textref h:C v:C \htext(8 26.5){\footnotesize$n-k-3$}

\move(36 24)\lvec(60 24)\lvec(60 20)\lvec(56 20)\lvec(56 16)\lvec(52
16)\lvec(52 12)\lvec(36 12)\lvec(36 24)

 \textref h:C v:C \htext(48 26.5){\footnotesize$\leq n+k-2$}
 \textref h:C v:C \htext(38 26){\footnotesize$\leftarrow$}
 \textref h:C v:C \htext(58 26){\footnotesize$\rightarrow$}
 \textref h:C v:C \htext(44 9.5){\footnotesize$\ge n-k+1$}

\move(80 24)\lvec(112 24)\lvec(112 20)\lvec(104 20)\lvec(104
16)\lvec(96 16)\lvec(96 12)\lvec(88 12)\lvec(88 8)\lvec(80
8)\lvec(80 24)

 \textref h:C v:C \htext(96 26.5){\footnotesize$\leq
2k$}
 \textref h:C
v:C \htext(83 26){\footnotesize$\longleftarrow$} \textref h:C v:C
\htext(109 26){\footnotesize$\longrightarrow$}

\linewd 0.1 \setgray 0

\move(36 24)\lvec(36 28)\move(60 24)\lvec(60 28)\move(0
24)\lvec(0 28)\move(16 24)\lvec(16 28) \move(80 24)\lvec(80
28)\move(112 24)\lvec(112 28)

\move(36 50)\lvec(36 46) \move(52 50)\lvec(52 46)
\move(36 12)\lvec(36 8) \move(52 12)\lvec(52 8)

\move(36
62)\lvec(36 66)\move(60 62)\lvec(60 66)\move(0 62)\lvec(0
66)\move(24 62)\lvec(24 66) \move(80 62)\lvec(80 66)\move(112
62)\lvec(112 66)

\textref h:C v:C \htext(12 4){\small$\tau'$}
 \textref h:C v:C \htext(48 4){\small$\lambda'$}
  \textref h:C v:C \htext(92 4){\small$\mu'$}

\textref h:C v:C \htext(12 42){\small$\tau$}
 \textref h:C v:C \htext(48 42){\small$\lambda$}
  \textref h:C v:C \htext(92 42){\small$\mu$}

\textref h:C v:C \htext(48 38){$\downarrow$}
  }
  \caption{The bijection $\varphi_A$ in Case $2$.}\label{fig-case2}
\end{figure}

\noindent Case $2$:
 $(\tau,\lambda,\mu)\in A_{n,k}$.
Removing the first two rows from $\tau$ and removing the first row from
$\mu$, we get $\tau'$ and $\mu'$, respectively. Let
$\lambda'=\lambda$. Then we have $(\tau',\lambda',\mu')\in A'_{n,k}$ and
\[
|\tau|+|\lambda|+|\mu|=2n-3+|\tau'|+|\lambda'|+|\mu'|.
\]
So we obtain a bijection $
 \varphi_A \colon A_{n,k}\rightarrow
 \{2n-3\}\times A'_{n,k}
$ as given by $ (\tau,\lambda,\mu )\mapsto (2n-3,
(\tau',\lambda',\mu')).$
 Figure \ref{fig-case2} gives an illustration of the correspondence.\\

\begin{figure}[h,t]
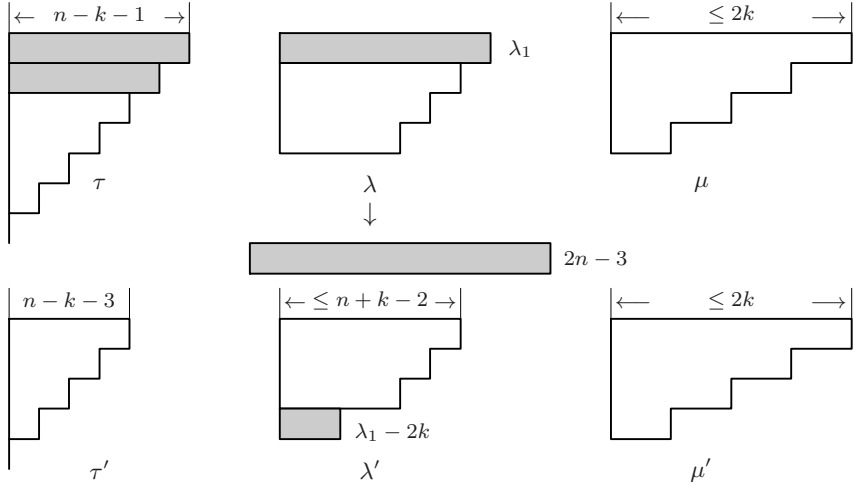

\centertexdraw{ \drawdim mm \linewd 0.25 \setgray 0

\move(0 54) \lvec(20 54) \lvec(20 58)\lvec(24 58)\lvec(24 62)\lvec(0
62)\lvec(0 54) \ifill f:0.8

\move(0 54) \lvec(20 54) \lvec(20 58)\lvec(24 58)\lvec(24 62)\lvec(0
62)\lvec(0 54) \move(0 58)\lvec(20 58) \move(0 54)\lvec(0 38)\lvec(4
38)\lvec(4 42)\lvec(8 42)\lvec(8 46)\lvec(12 46)\lvec(12 50)\lvec(16
50)\lvec(16 54)\move(0 38) \lvec(0 34)

\textref h:C v:C \htext(12 64.5){\footnotesize$n-k-1$}
 \textref h:C
v:C \htext(2 64){\footnotesize$\leftarrow$} \textref h:C v:C
\htext(22 64){\footnotesize$\rightarrow$}

\move(36 62)\lvec(64 62)\lvec(64 58)\lvec(36 58)\lvec(36 62)\ifill
f:0.8 \move(36 62)\lvec(64 62)\lvec(64 58)\lvec(36 58)\lvec(36 62)

\move(36 58)\lvec(60 58)\lvec(60 54)\lvec(56 54)\lvec(56 50)\lvec(52
50)\lvec(52 46)\lvec(36 46)\lvec(36 58)

 \textref h:C v:C \htext(68 60){\footnotesize$\lambda_1$}

\move(80 62)\lvec(112 62)\lvec(112 58)\lvec(104 58)\lvec(104
54)\lvec(96 54)\lvec(96 50)\lvec(88 50)\lvec(88 46)\lvec(80
46)\lvec(80 62)

 \textref h:C v:C \htext(96 64.5){\footnotesize$\leq
2k$}
 \textref h:C
v:C \htext(83 64){\footnotesize$\longleftarrow$} \textref h:C v:C
\htext(109 64){\footnotesize$\longrightarrow$}

\move(32 30)\lvec(72 30)\lvec(72 34)\lvec(32 34)\lvec(32 30) \ifill
f:0.8 \textref h:C v:C \htext(78 32){\footnotesize$2n-3$}

\move(32 30)\lvec(72 30)\lvec(72 34)\lvec(32 34)\lvec(32 30)

\move(0 24)\lvec(16 24)\lvec(16 20)\lvec(12 20 )\lvec(12 16)\lvec(8
16)\lvec(8 12)\lvec(4 12)\lvec(4 8)\lvec(0 8)\lvec(0 24)\move(0 8)
\lvec(0 4)

\textref h:C v:C \htext(8 26.5){\footnotesize$n-k-3$}

\move(36 24)\lvec(60 24)\lvec(60 20)\lvec(56 20)\lvec(56 16)\lvec(52
16)\lvec(52 12)\lvec(36 12)\lvec(36 24)

\move(36 12)\lvec(44 12)\lvec(44 8)\lvec(36 8)\lvec(36 12)\ifill
f:0.8 \move(36 12)\lvec(44 12)\lvec(44 8)\lvec(36 8)\lvec(36 12)

 \textref h:C v:C \htext(48
26.5){\footnotesize$\leq n+k-2$}
 \textref h:C
v:C \htext(38 26){\footnotesize$\leftarrow$} \textref h:C v:C
\htext(58 26){\footnotesize$\rightarrow$}
 \textref h:C v:C \htext(51 9){\footnotesize$\lambda_1-2k$}

\move(80 24)\lvec(112 24)\lvec(112 20)\lvec(104 20)\lvec(104
16)\lvec(96 16)\lvec(96 12)\lvec(88 12)\lvec(88 8)\lvec(80
8)\lvec(80 24)

 \textref h:C v:C \htext(96 26.5){\footnotesize$\leq
2k$}
 \textref h:C
v:C \htext(83 26){\footnotesize$\longleftarrow$} \textref h:C v:C
\htext(109 26){\footnotesize$\longrightarrow$}

\linewd 0.1 \setgray 0
\move(36 24)\lvec(36 28)\move(60 24)\lvec(60 28) \move(0 24)\lvec(0
28)\move(16 24)\lvec(16 28) \move(80 24)\lvec(80 28)\move(112
24)\lvec(112 28)

\move(0 62)\lvec(0 66)\move(24 62)\lvec(24 66) \move(80 62)\lvec(80
66)\move(112 62)\lvec(112 66)

\textref h:C v:C \htext(12 4){\small$\tau'$}
 \textref h:C v:C \htext(48 4){\small$\lambda'$}
  \textref h:C v:C
\htext(92 4){\small$\mu'$}

\textref h:C v:C \htext(12 42){\small$\tau$}
 \textref h:C v:C \htext(48 42){\small$\lambda$}
  \textref h:C v:C
\htext(92 42){\small$\mu$}

\textref h:C v:C \htext(48 38){$\downarrow$}
  }
  \caption{The bijection $\varphi_B$ in Case $3$.}\label{fig-case3}
\end{figure}

\noindent Case $3$: $(\tau,\lambda,\mu)\in B_{n,k}$. Removing the
first two rows from $\tau$, we get $\tau'$. Subtracting $2k$ from  the part $\lambda_1$ in
$\lambda$, we get a partition $\lambda'$. Let $\mu'=\mu$.
 Then we have  $(\tau',\lambda',\mu') \in  B'_{n,k}$ and
 \[
 |\tau|+|\lambda|+|\mu|=2n-3+|\tau'|+|\lambda'|+|\mu'|.
\]
 Thus we obtain a bijection $
 \varphi_B\colon B_{n,k}\rightarrow
 \{2n-3\}\times B'_{n,k}
$
defined by
$
 (\tau,\lambda,\mu )\mapsto (2n-3,  (\tau',\lambda',\mu')).
$
See Figure \ref{fig-case3} for an illustration.\\

\begin{figure}[h,t]
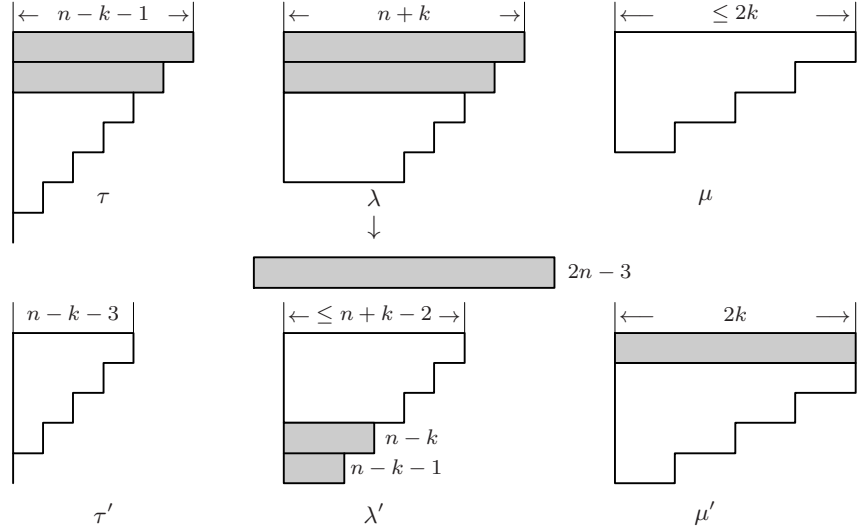

\centertexdraw{ \drawdim mm \linewd 0.25 \setgray 0

\move(0 56) \lvec(20 56) \lvec(20 60)\lvec(24 60)\lvec(24 64)\lvec(0
64)\lvec(0 56) \ifill f:0.8

\move(0 56) \lvec(20 56) \lvec(20 60)\lvec(24 60)\lvec(24 64)\lvec(0
64)\lvec(0 56)  \move(0 60)\lvec(20 60) \move(0 56)\lvec(0
40)\lvec(4 40)\lvec(4 44)\lvec(8 44)\lvec(8 48)\lvec(12 48)\lvec(12
52)\lvec(16 52)\lvec(16 56)\move(0 40) \lvec(0 36)

\textref h:C v:C \htext(12 66.5){\footnotesize$n-k-1$} \textref h:C
v:C \htext(2 66){\footnotesize$\leftarrow$} \textref h:C v:C
\htext(22 66){\footnotesize$\rightarrow$}

\move(36 56)\lvec(60 56)\lvec(60 52)\lvec(56 52)\lvec(56 48)\lvec(52
48)\lvec(52 44)\lvec(36 44)\lvec(36 56)

\move(36 64)\lvec(68 64)\lvec(68 60)\lvec(64 60) \lvec(64 56)
\lvec(36 56) \lvec(36 64)  \ifill f:0.8

\move(36 60)\lvec(64 60) \move(36 64)\lvec(68 64)\lvec(68
60)\lvec(64 60) \lvec(64 56) \lvec(36 56) \lvec(36 64)

\textref h:C v:C \htext(52 66.5){\footnotesize$ n+k$} \textref h:C
v:C \htext(38 66){\footnotesize$\leftarrow$} \textref h:C v:C
\htext(66 66){\footnotesize$\rightarrow$}

\move(80 64)\lvec(112 64)\lvec(112 60)\lvec(104 60)\lvec(104
56)\lvec(96 56)\lvec(96 52)\lvec(88 52)\lvec(88 48)\lvec(80
48)\lvec(80 64)

\textref h:C v:C \htext(96 66.5){\footnotesize$\leq 2k$} \textref
h:C v:C \htext(83 66){\footnotesize$\longleftarrow$} \textref h:C
v:C \htext(109 66){\footnotesize$\longrightarrow$}

\move(32 30)\lvec(72 30)\lvec(72 34)\lvec(32 34)\lvec(32 30) \ifill
f:0.8 \textref h:C v:C \htext(78 32){\footnotesize$2n-3$}

\move(32 30)\lvec(72 30)\lvec(72 34)\lvec(32 34)\lvec(32 30)

\move(0 24)\lvec(16 24)\lvec(16 20)\lvec(12 20 )\lvec(12 16)\lvec(8
16)\lvec(8 12)\lvec(4 12)\lvec(4 8)\lvec(0 8)\lvec(0 24)\move(0 8)
\lvec(0 4)

\textref h:C v:C \htext(8 26.5){\footnotesize$n-k-3$}

\move(36 24)\lvec(60 24)\lvec(60 20)\lvec(56 20)\lvec(56 16)\lvec(52
16)\lvec(52 12)\lvec(36 12)\lvec(36 24)

\move(36 12)\lvec(48 12)\lvec(48 8)\lvec(44 8)\lvec(44 4)\lvec(36
4)\lvec(36 12)\ifill f:0.8 \move(36 12)\lvec(48 12)\lvec(48
8)\lvec(44 8)\lvec(44 4)\lvec(36 4)\lvec(36 12) \move(36 8)\lvec(44
8)

\textref h:C v:C \htext(48 26.5){\footnotesize$\leq n+k-2$} \textref
h:C v:C \htext(38 26){\footnotesize$\leftarrow$} \textref h:C v:C
\htext(58 26){\footnotesize$\rightarrow$} \textref h:C v:C \htext(53
10){\footnotesize$ n-k$} \textref h:C v:C \htext(51
6){\footnotesize$ n-k-1$}

\move(80 20)\lvec(112 20)\lvec(112 16)\lvec(104 16)\lvec(104
12)\lvec(96 12)\lvec(96 8)\lvec(88 8)\lvec(88 4)\lvec(80 4)\lvec(80
20)

\move(80 24)\lvec(112 24)\lvec(112 20)\lvec(80 20)\lvec(80 24)\ifill
f:0.8 \move(80 24)\lvec(112 24)\lvec(112 20)\lvec(80 20)\lvec(80 24)

 \textref h:C v:C \htext(96 26.5){\footnotesize$2k$}
 \textref h:C
v:C \htext(83 26){\footnotesize$\longleftarrow$} \textref h:C v:C
\htext(109 26){\footnotesize$\longrightarrow$}

\linewd 0.1 \setgray 0 \move(0 64)\lvec(0 68)\move(24 64)\lvec(24
68) \move(36 64)\lvec(36 68)\move(68 64)\lvec(68 68) \move(80
64)\lvec(80 68)\move(112 64)\lvec(112 68) \move(0 24)\lvec(0
28)\move(16 24)\lvec(16 28) \move(36 24)\lvec(36 28) \move(60
24)\lvec(60 28) \move(80 24)\lvec(80 28)\move(112 24)\lvec(112 28)

\textref h:C v:C \htext(12 0){\small$\tau'$}
 \textref h:C v:C \htext(48 0){\small$\lambda'$}
  \textref h:C v:C
\htext(92 0){\small$\mu'$}

\textref h:C v:C \htext(12 42){\small$\tau$}
 \textref h:C v:C \htext(48 42){\small$\lambda$}
  \textref h:C v:C
\htext(92 42){\small$\mu$}

\textref h:C v:C \htext(48 38){$\downarrow$}
  }
  \caption{The bijection $\varphi_C$ in Case $4$.}\label{fig-case4}
\end{figure}

\noindent Case $4$: $(\tau,\lambda,\mu)\in C_{n,k}$. Removing
first two rows from $\tau$ we get $\tau'$.  Subtracting  $2k$ from the
parts $n+k-1$ and $n+k$ in
$\lambda$, we get a partition $\lambda'$.
Adding $2k$ to $\mu$ as a new part, we get $\mu'$.  Then we have
$(\tau',\lambda',\mu') \in
 C'_{n,k} $ and
\[
|\tau|+|\lambda|+|\mu|=2n-3+|\tau'|+|\lambda'|+|\mu'|.
\]
 Thus we obtain  a bijection $
 \varphi_C\colon C_{n,k}\rightarrow
 \{2n-3\}\times C'_{n,k}
$
as given by
$
 (\tau,\lambda,\mu )\mapsto (2n-3,  (\tau',\lambda',\mu')).
$ This case is illustrated in Figure \ref{fig-case4}.

\begin{figure}[h,t]
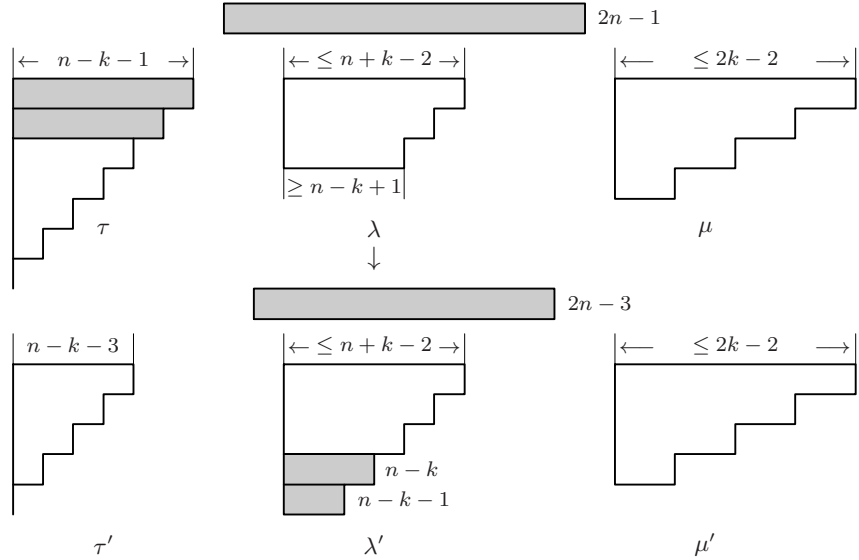

\centertexdraw{ \drawdim mm \linewd 0.25 \setgray 0

 \move(28
68)\lvec(76 68)\lvec(76 72)\lvec(28 72)\lvec(28 68) \ifill f:0.8
\textref h:C v:C \htext(82 70){\footnotesize$2n-1$}

 \move(28
68)\lvec(76 68)\lvec(76 72)\lvec(28 72)\lvec(28 68)

\move(0 54) \lvec(20 54) \lvec(20 58)\lvec(24 58)\lvec(24 62)\lvec(0
62)\lvec(0 54) \ifill f:0.8

\move(0 54) \lvec(20 54) \lvec(20 58)\lvec(24 58)\lvec(24 62)\lvec(0
62)\lvec(0 54) \move(0 58)\lvec(20 58) \move(0 54)\lvec(0 38)\lvec(4
38)\lvec(4 42)\lvec(8 42)\lvec(8 46)\lvec(12 46)\lvec(12 50)\lvec(16
50)\lvec(16 54)\move(0 38) \lvec(0 34)

\textref h:C v:C \htext(12 64.5){\footnotesize$n-k-1$}
 \textref h:C
v:C \htext(2 64){\footnotesize$\leftarrow$} \textref h:C v:C
\htext(22 64){\footnotesize$\rightarrow$}

\move(36 62)\lvec(60 62)\lvec(60 58)\lvec(56 58)\lvec(56 54)\lvec(52
54)\lvec(52 50)\lvec(36 50)\lvec(36 62)

 \textref h:C v:C \htext(38 64){\footnotesize$\leftarrow$}
 \textref h:C v:C \htext(48 64.5){\footnotesize$\leq n+k-2$}
 \textref h:C v:C \htext(58 64){\footnotesize$\rightarrow$}
 \textref h:C v:C \htext(44 47.5){\footnotesize$\ge n-k+1$}

\move(80 62)\lvec(112 62)\lvec(112 58)\lvec(104 58)\lvec(104
54)\lvec(96 54)\lvec(96 50)\lvec(88 50)\lvec(88 46)\lvec(80
46)\lvec(80 62)

 \textref h:C v:C \htext(96 64.5){\footnotesize$\leq
2k-2$}
 \textref h:C
v:C \htext(83 64){\footnotesize$\longleftarrow$} \textref h:C v:C
\htext(109 64){\footnotesize$\longrightarrow$}

\move(32 30)\lvec(72 30)\lvec(72 34)\lvec(32 34)\lvec(32 30) \ifill
f:0.8 \textref h:C v:C \htext(78 32){\footnotesize$2n-3$}

\move(32 30)\lvec(72 30)\lvec(72 34)\lvec(32 34)\lvec(32 30)

\move(0 24)\lvec(16 24)\lvec(16 20)\lvec(12 20 )\lvec(12 16)\lvec(8
16)\lvec(8 12)\lvec(4 12)\lvec(4 8)\lvec(0 8)\lvec(0 24)\move(0 8)
\lvec(0 4)

\textref h:C v:C \htext(8 26.5){\footnotesize$n-k-3$}

\move(36 24)\lvec(60 24)\lvec(60 20)\lvec(56 20)\lvec(56 16)\lvec(52
16)\lvec(52 12)\lvec(36 12)\lvec(36 24)

\move(36 12)\lvec(48 12)\lvec(48 8)\lvec(44 8)\lvec(44 4)\lvec(36
4)\lvec(36 12)\ifill f:0.8 \move(36 12)\lvec(48 12)\lvec(48
8)\lvec(44 8)\lvec(44 4)\lvec(36 4)\lvec(36 12) \move(36 8)\lvec(44
8)

 \textref h:C v:C \htext(48 26.5){\footnotesize$\leq n+k-2$}
 \textref h:C v:C \htext(38 26){\footnotesize$\leftarrow$}
 \textref h:C v:C \htext(58 26){\footnotesize$\rightarrow$}
 \textref h:C v:C \htext(52 6){\footnotesize$n-k-1$}
 \textref h:C v:C \htext(53 10){\footnotesize$n-k$}

\move(80 24)\lvec(112 24)\lvec(112 20)\lvec(104 20)\lvec(104
16)\lvec(96 16)\lvec(96 12)\lvec(88 12)\lvec(88 8)\lvec(80
8)\lvec(80 24)

 \textref h:C v:C \htext(96 26.5){\footnotesize$\leq
2k-2$}
 \textref h:C
v:C \htext(83 26){\footnotesize$\longleftarrow$} \textref h:C v:C
\htext(109 26){\footnotesize$\longrightarrow$}

\linewd 0.1 \setgray 0
\move(36 24)\lvec(36 28)\move(60 24)\lvec(60 28)\move(0 24)\lvec(0
28)\move(16 24)\lvec(16 28) \move(80 24)\lvec(80 28)\move(112
24)\lvec(112 28)

\move(36 50)\lvec(36 46) \move(52 50)\lvec(52 46) \move(36
62)\lvec(36 66)\move(60 62)\lvec(60 66) \move(0 62)\lvec(0
66)\move(24 62)\lvec(24 66) \move(80 62)\lvec(80 66)\move(112
62)\lvec(112 66)

\textref h:C v:C \htext(12 0){\small$\tau'$}
 \textref h:C v:C \htext(48 0){\small$\lambda'$}
  \textref h:C v:C
\htext(92 0){\small$\mu'$}

\textref h:C v:C \htext(12 42){\small$\tau$}
 \textref h:C v:C \htext(48 42){\small$\lambda$}
  \textref h:C v:C
\htext(92 42){\small$\mu$}

\textref h:C v:C \htext(48 38){$\downarrow$}
  }
  \caption{The bijection $\varphi_D$ on $\{2n-1\}\times P_{n-1,k-1}$.}\label{fig-caseD}
\end{figure}

Now we consider the quadruples
\[
(2n-1,(\tau,\lambda,\mu))\in \{2n-1\}\times
P_{n-1,k-1}.
\]
For any $(\tau,\lambda,\mu )\in P_{n-1,k-1}$,  remove
the first two rows of $\tau$ and  add two parts $n-k$ and $n-k-1$ to
$\lambda$ to get $\tau'$ and $\lambda'$. Let $\mu'=\mu$.
Then we see that $(\tau',\lambda',\mu')\in D_{n,k}$ and
\[
2n-1+|\tau|+|\lambda|+|\mu|=2n-3+|\tau'|+|\lambda'|+|\mu'|.
\]
Thus we obtain a bijection
\[
 \varphi_D \colon \{2n-1\}\times P_{n-1,k-1}\rightarrow
 \{2n-3\}\times D_{n,k}
\] as given by $ (2n-1,(\tau,\lambda,\mu ))\mapsto (2n-3,
(\tau',\lambda',\mu')).
$ This case is illustrated by Figure \ref{fig-caseD}.\\

The proof is complete by combining the bijections $\varphi_A$, $\varphi_B$, $\varphi_C$ and $\varphi_D$. \qed

The above theorem gives the bijections $\phi_{n,k}$ for $0\le k \le n-2$.
In the   following theorem we consider the special cases $k=n-1$ and $k=n$.

\begin{theo}\label{thm-involution}
For $n\geq 2$ and for $k=n-1$ or $n$, there is an
involution
\begin{equation}\label{involution-q12}
 I_{n,k}\colon  P_{n,k} \cup \{2n-1\}\times P_{n-1,k-1}
              \rightarrow P_{n-1,k-1}.
 \end{equation}
\end{theo}

\pf We only give the description  of the involution
 $I_{n,n}$ since  $I_{n,n-1}$ can
be constructed in the same manner.

\noindent Case $1.$  For $(\emptyset,\lambda,\mu )\in P_{n,n}$,
if the first part of $\lambda$
is $2n$, then move it to $\mu$. Conversely, if $\mu$ contains a part $2n$
but $\lambda$ does not,
then move this part from $\mu$ back to $\lambda$.

\noindent Case $2.$  For  $(\emptyset, \lambda,\mu)\in P_{n,n}$ with
$\lambda_1=2n-1$ and $\mu_1<2n$, remove
 the first part $2n-1$ of $\lambda$ to get $\lambda'$,
and set
\[  I_{n,n}(\emptyset, \lambda,\mu)=(2n-1,(\emptyset,\lambda',\mu),\]
 which belongs to
$\{2n-1\}\times P_{n-1,n-1}$. Conversely, for
\[(2n-1,(\emptyset, \lambda,\mu))\in
\{2n-1\}\times P_{n-1,n-1},\]
add  a part $2n-1$ to $\lambda$, we get
$\lambda'$ and set
\[ I_{n,n}(2n-1, (\emptyset, \lambda, \mu))= (\emptyset, \lambda',\mu),\]
 which belongs to
$P_{n,n}$.

\noindent Case $3.$  It can be seen that the set of triples
 $(\emptyset, \lambda,\mu)\in P_{n,n}$ with
 $\lambda_1<2n-1$ and $\mu_1<2n$  is exactly $P_{n-1,n-1}$. So we set
 $P_{n-1, n-1}$ to be the invariant set of the involution.

Thus we obtain an  involution on $P_{n,n} \cup \{2n-1\}\times
P_{n-1,n-1}$ with the invariant set
$P_{n-1,n-1}$.  \qed

Define a weight function $w$ on $P_{n,k}$, $\{2n-1\}\times
P_{n-1,k}$ and $\{2n-3\}\times P_{n-2,k}$ as given by
 \begin{eqnarray*}
  w(\tau,\lambda,\mu)&=&(-1)^{\ell(\lambda)}\,
          q ^{\mid\tau\mid+\mid\lambda\mid+\mid\mu\mid},\\[5pt]
  w(2n-1,(\tau,\lambda,\mu))&=& q^{2n-1}\, (-1)^{\ell(\lambda)} \,
          q ^{\mid\tau\mid+\mid\lambda\mid+\mid\mu\mid},\\[5pt]
  w(2n-3,(\tau,\lambda,\mu))&=& q^{2n-3}\, (-1)^{\ell(\lambda)} \,
             q ^{\mid\tau\mid+\mid\lambda\mid+\mid\mu\mid}.
 \end{eqnarray*}
 One sees that the bijections and involutions in Theorems \ref{thm-bijection} and \ref{thm-involution}  are
 weight
preserving.  Hence we get the following recurrence relation for
\[
F_n(q)=\sum_{k\ge 0} F_{n,k}.
\]

\begin{coro}
For $n\geq 2$,  we have
\begin{equation}\label{rec-q12-final'}
F_n(q)+(q^{2n-1}-1)F_{n-1}(q)-q^{2n-3}F_{n-2}(q)=0.
\end{equation}
\end{coro}

It is easy to verify that the right hand side of (\ref{eq-q12a}), namely, the sum
\begin{equation}\label{eq-q12-final'}
(-1)^nq^{n^2}\sum_{j=-n}^{n}(-1)^jq^{-j^2},
\end{equation}
also satisfies the recurrence relation (\ref{rec-q12-final'}).
Taking the initial values into consideration, we are led to the identity
of Andrews.

\vspace{.2cm} \noindent{\bf Acknowledgments.}
This work was supported by the 973 Project, the PCSIRT Project of the Ministry of
Education, and the National Science Foundation of China.


\begin{thebibliography}{99}

\bibitem{Andrews1998}
G.E.~Andrews, The Theory of Partitions, Cambridge University Press,
 Cambridge, 1998.


\bibitem{Andrews2010}
G.E.~Andrews, Parity in partition identities, Ramanujan J. 23 (2010)
45--90.

\bibitem{Chen2005}
W.Y.C.~Chen, Q.-H.~Hou  and Lisa.H.~Sun, The method of combinatorial telescoping, J. Combin. Theory, Ser. A 118 (2011) 899--907.

\bibitem{FP}
D. Feldman  and J. Propp,
 Producing new bijections from old, Adv. Math.
113 (1995) 1--44.


\bibitem{Gasper}
G.~Gasper and M.~Rahman, Basic Hypergeometric Series, Encyclopedia
of Mathematics and Its Applications, Vol. 35, Cambridge University
Press, Cambridge, 1990.


\bibitem{Knuth1994}
R.~Graham, D.~Knuth, O.~Patashnik, Concrete Mathematics, 2nd
Ed., Addison-Wesley, Reading, MA, 1994.

\bibitem{MacMahon}
P.A.~MacMahon, Combinatory Analysis, Cambridge University Press,
Cambridge, (1916).

\bibitem{Pak2006}
I.~Pak, Partition bijectionis, a survey, Ramanujan J. 12 (2006) 50--57.

\bibitem{A=B1996}
M.~Petkov\v{s}ek, H.S.~Wilf, and D.~Zeilber, $A=B$,  A.K.~Peters,
Wellesley, MA, 1996.



\bibitem{Watson}
G.N.~Watson, A new proof of the Rogers-Ramanujan identities, J.
London Math. Soc. 4 (1929) 4--9.



\bibitem{Syl}
J.J.~Sylvester, A constructive theory of partitions, arranged in
three acts, an interact, and an exodion, Amer. J. Math. 5 (1882)
251--330.



\bibitem{Zeilberger1991}
D.~Zeilberger, The method of creative telescoping, J. Symbolic
Comput. 11 (1991) 195--204.
\end{thebibliography}
\end{document}